\documentclass[12pt]{article}
\usepackage{amssymb}
\usepackage{latexsym}

\newtheorem{theorem}{Theorem}[section]
\newtheorem{lemma}[theorem]{Lemma}
\newtheorem{definition}[theorem]{Definition}
\newtheorem{proposition}[theorem]{Proposition}
\newtheorem{corollary}[theorem]{Corollary}
\newtheorem{notation}[theorem]{Notation}

\newenvironment{proof}{\paragraph{\it Proof.}}{$\square$\vskip0.4cm}
\newenvironment{remark}{\paragraph{\it Remark.}}{\vskip0.4cm}
\newenvironment{example}{\paragraph{\it Example.}}{\vskip0.4cm}

\newcommand{\M}{{\cal M}}
\newcommand{\N}{{\cal N}}
\newcommand{\C}{{\mathbf C}}
\newcommand{\R}{{\mathbf R}}
\newcommand{\Proj}{{\mathbf P}}
\newcommand{\Z}{{\mathbf Z}}
\newcommand{\cM}{{{\bar\M}}}
\newcommand{\calO}{{\cal O}}
\newcommand{\calL}{{\cal L}}
\newcommand{\Pic}{\mathop{\rm Pic}\nolimits}
\newcommand{\coker}{\mathop{\rm coker}\nolimits}
\newcommand{\sing}{\mathop{\rm sing}\nolimits}

\begin{document}

\title{Compactification of moduli of Higgs bundles}
\author{Tam\'as Hausel\\ DPMMS, University of Cambridge}
\maketitle

\begin{abstract} 

In this paper we consider a canonical compactification of  
$\M$, the moduli space of stable Higgs bundles with 
fixed determinant of odd degree
over a Riemann surface $\Sigma$, producing a projective variety
$\cM=\M\cup Z$. We give a detailed study of the spaces $\cM$,
$Z$ and $\M$. In doing so we reprove some assertions of Laumon and Thaddeus on
the nilpotent cone.   
\end{abstract}

\section{Introduction}

Magnetic monopoles, the solutions of Bogomolny equations of mathematical 
physics, can be interpreted as solutions of the self-dual Yang-Mills 
equations on $\R^4$ which are translation invariant in
one direction.
Motivated by this interpretation, 
Hitchin in {\cite{hitchin1}} addressed the problem of finding  
solutions
to the $SU(2)$ self-dual Yang-Mills  equations on $\R^4$, 
which are translation invariant in two directions. Although such
solutions of finite energy do not exist, 
due to  the conformal invariance of the equations,
it was possible to find
solutions of the corresponding $SU(2)$ self-duality 
equations over a Riemann surface
$\Sigma$. In the same paper Hitchin gave an extensive description of
the space of these solutions. 

One important result shows how to 
assign in a certain one to one manner an algebro-geometric object
to a solution of Hitchin's self-duality equations. This algebro-geometric
object is called a stable Higgs bundle, which consists of a pair of
a rank 2 holomorphic vector bundle $E$ on  $\Sigma$ and a section 
$\Phi \in H^0(\Sigma,End_0(E)\otimes K_{\Sigma})$. The latter is 
called the Higgs field, after the analogous object in the monopole case.  

In \cite{hitchin1}, in these algebro-geometric terms, Hitchin investigates 
the moduli space $\M$ of stable Higgs bundles with fixed determinant of
degree 1. This notion and the corresponding moduli space has become 
important from a purely algebro-geometric point of view, too. The main reason
is that the cotangent bundle of $\N$, the moduli space of stable 
rank 2 vector bundles with fixed determinant of odd degree, which
is a well researched object in the algebraic geometry of vector bundles,
sits inside
$\M$ as an open dense subset. Namely, $(T^*_\N)_E$ is canonically isomorphic
to $H^0(\Sigma, End_0(E)\otimes K_\Sigma)$ thus the points of $T^*_\N$
are Higgs bundles. 

Among other results Hitchin proved that $\M$ is a non-compact complete
hyperk\"ahler manifold. Defined as above, in purely 
algebro-geometric terms, it 
was not surprising that
$\M$ turned out to be quasi-projective as Nitsure
has shown in \cite{nitsure}.

The main aim of this paper is to investigate a canonical compactification of
$\M$: among other things we show that the compactification is projective,
calculate its Picard group, and calculate the Poincar\'e polynomial for the
cohomology.

In this paper we use a simple method to compactify non-compact
K\"ahler manifolds with a nice proper Hamiltonian $S^1$ action via
Lerman's construction of symplectic cutting \cite{larman}. We use
this method to compactify $\M$.
Our approach is symplectic in nature and eventually
produces some fundamental results about the spaces occurring, using existing
techniques from the theory of symplectic quotients. 

We show that the compactification described in this paper 
is a good example of  Yau's problem of
finding a complete Ricci flat metric on the complement of 
a nef anticanonical divisor
in a projective variety. 

Many of the results of this paper can be easily generalized
to  other Higgs bundle moduli spaces, which have been extensively investigated
(see e.g. \cite{nitsure} and \cite{simpson}). As a matter of fact 
Simpson gave a 
definition of a similar compactificitation for these more general Higgs bundle
moduli spaces in Theorems 11.2 and 11.1 of \cite{simpson2} and 
in Proposition 17 of \cite{simpson3}, without investigating it in detail. 
For example, the projectiveness of the compactification 
is not clear from these definitions. One novelty of our paper is the
proof of the projectiveness of the compactification in our case.

Since the compactification method used in this paper is fairly 
general it is possible
to apply it to other K\"ahler manifolds with the above properties. It could
be interesting for instance to see how this method works for
the toric hyperk\"ahler manifolds of Goto \cite{goto} and Bielawski and Dancer
\cite{bielawski-dancer}.

Finally, as a conclusion, we note that the compactification of this paper 
solves one half of the problem
of compactifying the moduli space $\M$, namely the `outer' half, i.e. shows
what the resulting spaces look like; while the other half of the problem
the `inner' part, i.e. how this fits into the moduli space description
of $\M$, is treated in the recent paper of Schmitt \cite{schmitt}. Schmitt's
approach is algebro-geometric in nature, and concerns mainly the construction
of the right notion for moduli to produce $\cM$, thus complements the present
paper.    
The relation between the two approaches deserves further investigation.  

\vskip0.4cm 

\paragraph{\bf Acknowledgements.} First of all I would like to thank
my supervisor Nigel Hitchin for fruitful supervisions. 
The compatibility with Yau's problem was
suggested by Michael Atiyah, while  Lerman's symplectic cutting was suggested 
by Michael Thaddeus. With both of them I had very inspiring conversations. 
I also thank Bal\'azs Szendr\H oi and the referee for helpful comments. 
Finally, I thank Trinity College, Cambridge for financial support.

\section{Statement of results}

In this section we describe the structure of the paper and list 
the results.

In Section~\ref{rizsa}  we collect the existing 
results about $\C^*$ actions on K\"ahler
manifolds and 
subsequently on K\"ahler quotients from the literature. 
We explain a general 
method of compactifying K\"ahler manifolds  with a nice, proper,
Hamiltonian $S^1$ action. The rest of the paper follows the structure
of Section~\ref{rizsa}.

In Section~\ref{higgs} we define the basic notions and restate some 
results of Hitchin about $\M$. Here we learn that the results of 
Section~\ref{rizsa} apply to $\M$. We describe here our toy example 
$\M_{toy}$, the moduli space of parabolic Higgs bundles on $\Proj^1_4$, which
serves as an example throughout the paper. 

In Section~\ref{nilpotent} (following ideas of Subsection~\ref{rizsa1}) 
we describe the nilpotent cone after Thaddeus
\cite{thaddeus1} 
and show that it coincides with the downward Morse
flow (Theorem~\ref{morse}). We reprove 
Laumon's theorem in our case, that the nilpotent cone is Lagrangian 
(Corollary~\ref{lagrange}). 

In Section~\ref{kaehler} we describe $Z$, the highest level K\"ahler quotient
of $\M$, while in \ref{kompakt} we analyse $\cM=\M\cup Z$. 
Here we follow the approaches of Subsection~\ref{rizsa2} and
Subsection~\ref{rizsa3}, respectively. 
Among others, we prove the following
statements:

\begin{itemize}

\item

$\cM$ is a compactification of $\M$, the moduli space of stable Higgs bundles
with fixed determinant and degree $1$ (Theorem~\ref{compactification}).

\item

$Z$ is a symplectic quotient of $\cal M$ by the circle action
$(E,\Phi)\mapsto (E,e^{i\theta}\Phi)$. $\bar{\cal M}$ is a symplectic
quotient  of ${\cal M}\times{\mathbf C}$ with respect to the circle action,
which is the usual one on $\cal M$ and multiplication on $\mathbf C$.

\item
While $\M$ is a smooth manifold, $Z$ is an orbifold, with only
${\mathbf Z}_2$
singularities corresponding to the fixed point set of the map
$(E,\Phi)\mapsto (E,-\Phi)$ on $\M$ (Theorem~\ref{Zorbi}), 
while similarly $\cM$ is
an orbifold with only ${\mathbf Z}_2$ singularities, and the singular locus of
$\cM$ coincides with that of $Z$ (Theorem~\ref{Morbi}).

\item

The Hitchin map $$\chi : \M \rightarrow {\mathbf C}^{3g-3}$$ extends to
a map $$\bar{\chi}:\cM\rightarrow {\mathbf P}^{3g-3}$$ which when
restricted
to $Z$ gives a map $$\bar{\chi}: Z \rightarrow {\mathbf P}^{3g-4}$$ whose generic
fibre is a Kummer variety corresponding to the Prym variety of the generic
fibre of the Hitchin map (Theorem~\ref{Zchi}, Theorem~\ref{Mchi}).
\item

$\cM$ is a projective variety (Theorem~\ref{Mproj}), with divisor $Z$ such that
$$(3g-2)Z=-K_{\cM},$$
the anticanonical divisor of $\cM$ (Corollary~\ref{vonat}).

\item

Moreover, $Z$ itself is a projective variety (Theorem~\ref{proj})
with an inherited holomorphic contact structure with
contact line bundle $L_Z$ (Theorem~\ref{contact})
and a one-parameter family
of K\"ahler forms $\omega_t(Z)$ (Theorem~\ref{csalad}). The Picard group
of $Z$ is described in Corollary~\ref{picard}. Moreover, the normal bundle
of $Z$ in $\cM$ is $L_Z$ which is nef by Corollary~\ref{triv}.

\item

Furthermore, $\cM$ has a one-parameter family of
K\"ahler forms $\omega_t(\cM)$,
which when restricted $Z$ gives the above $\omega_t(Z)$. 

\item

$Z$ is birationally equivalent to $P(T^*_\N)$ the projectivized
cotangent bundle of the moduli space of rank 2 stable bundles with fixed
determinant and odd degree (Corollary~\ref{Zjeno}). 
$\bar{\M}$ is birationally equivalent to
$P(T^*_\N\oplus {\cal O}_\N)$, the canonical compactification of
$T^*_{\N}$ (Corollary~\ref{Mjeno}).

\item

We calculate certain sheaf cohomology groups in Corollary~\ref{H0Z} and
Corollary~\ref{H1Z} and interpret some of these results as the equality
of certain infinitesimal deformation spaces. 

\item

The Poincar\'e polynomial of $Z$ is described in Corollary~\ref{PZ}, the
Poincar\'e polynomial of $\cM$ is described in Theorem~\ref{PM}. 

\item

We finish Section~\ref{kompakt} by showing an interesting isomorphism between
two vector spaces: one contains information about the 
intersection of the components
of the nilpotent cone, the other says something 
about the contact line bundle $L_Z$ on $Z$.  
\end{itemize}

\section{Compactification by symplectic cutting}
\label{rizsa}

In this section we collect the results from the literature concerning
$\C^*$ actions on K\"ahler manifolds. At the same time we sketch the
structure of the rest of the paper.

\subsection{Stratifications}
\label{rizsa1}

Suppose that we are given a K\"ahler manifold $(M,I,\omega)$ with
complex structure $I$ and K\"ahler form $\omega$. Suppose also that
$\C^*$ acts on $M$
biholomorphically with respect to $I$ and such that the K\"ahler form is 
invariant under the induced action of $S^1\subset \C^*$. Suppose furthermore
that this latter action is Hamiltonian with proper moment map 
$\mu:M\rightarrow \R$, with finitely many critical points and 
$0$ being the absolute minimum of $\mu$.    
Let $\{N_\alpha\}_{\alpha\in A}$ be the set of the components of the fixed 
point
set of the $\C^*$ action. 

We list some results of \cite{kirwan} extended to our case. Namely, Kirwan's
results are stated for compact K\"ahler manifolds, but one can always modify
the proofs for non-compact manifolds as above (cf. Chapter 9 in 
\cite{kirwan}). 
  
There exist two stratifications in such a situation. The first one
is called the {\em Morse stratification} and can be defined as follows. 
The stratum $S^M_\alpha$ is the set of points of $\M$ whose path of 
steepest descent for the Morse function $\mu$ and
the K\"ahler metric have limit points in $N_\alpha$. One can also define 
the sets $T^M_\alpha$ as the points of $\M$ whose path of steepest descent
for the Morse function $-\mu$ and the K\"ahler metric have limit points in 
$N_\alpha$. $S^M_\alpha$ gives a stratification even in the 
non-compact case, however the set $\bigcup_\alpha
T^M_\alpha$is not the whole space but a 
deformation retract of it. The set $\bigcup_\alpha T^M_\alpha$ is called the 
{\em downward Morse flow}. 

The other stratification is the {\em Bialynicki-Birula stratification}, where
the stratum $S^B_\alpha$ is the set of points $p\in M$ for which 
$\lim_{t\rightarrow 0}tp\in N_\alpha$. 
Similarly, as above, we can define
$T^B_\alpha$ as the points $p\in M$ for which 
$\lim_{t\rightarrow \infty}tp\in N_\alpha$. 

One of Kirwan's important results in \cite{kirwan} Theorem $6.16$ 
asserts that the 
stratifications $S^M_\alpha$ and $S^B_\alpha$ coincide, and similarly 
$T^M_\alpha=T^B_\alpha=T_\alpha$. This result is important
because it shows that the strata $S_\alpha=S^M_\alpha=S^B_\alpha$ 
of the stratifications are total spaces
of affine bundles (so-called $\beta$-fibrations) on $N_\alpha$
(this follows 
from the Bialynicki-Birula picture) and moreover this stratification
is responsible for the topology of the space $M$ (this follows 
from the Morse 
picture). Thus we have the following theorem (cf. 
Theorem $4.1$ of \cite{bialynicki} and also Theorem $1.12$ of 
\cite{thaddeus3}): 

\begin{theorem} $S_\alpha$ and $T_\alpha$ are complex submanifolds 
of $M$. They are isomorphic to total spaces of some $\beta$-fibrations
over $N_\alpha$, such that the normal bundle of $N_\alpha$ in these 
$\beta$-fibrations
are $E^+_\alpha$ and $E^-_\alpha$, respectively, 
where $E^+_\alpha$ is the positive and $E^-_\alpha$ is the negative 
subbundle of $T_M\mid_{N_\alpha}$ with respect to the $S^1$ action. 

Moreover, the downward Morse flow $\bigcup_\alpha T_\alpha$ is a deformation
retract of $M$.
\label{stratification}
\end{theorem}

Recall that a $\beta$-fibration in our case 
is a fibration $E\rightarrow B^n$ with a $\C^*$ action on the 
total space which is locally like $\C^n \times V$, where $V$ is the $\C^*$ 
module $\beta:\C^*\rightarrow GL(V)$. Note that such a fibration is not
a vector bundle in general, but it is if $\beta$ is the sum of isomorphic,
one-dimensional non-trivial $\C^*$ modules.

\subsection{K\"ahler quotients} 
\label{rizsa2}

We define an action to be {\em semi-free} if the stabilizer of any point
is finite or the whole group itself. 
 
Whenever we are given a Hamiltonian,
proper, semi-free $S^1$ action on a K\"ahler manifold, we can form the K\"ahler quotients
$Q_t=\mu^{-1}(t)/S^1$, which are compact K\"ahler orbifolds at a regular value
$t$ of $\mu$. 

If this $S^1$ action is induced from an action of $\C^*$ on $M$ as  
above, then we can relate the K\"ahler quotients to the 
quotients $M/\C^*$ as follows. First we define 
$M^{min}_t\subset M$ as the set of points in $M$ whose $\C^*$ orbit 
intersects $\mu^{-1}(t)$. Now a theorem of Kirwan states (see Theorem 7.4
in \cite{kirwan}) that it is possible to define a complex structure on the
orbit space $M^{min}_t/\C^*$, 
and she also proves that this space is homeomorphic to $Q_t$, 
defining the complex structure for the K\"ahler quotient $Q_t$. (Here again
we used the results of Kirwan for non-compact manifolds, but as above, these 
results can be easily modified for our situation.) It now
simply follows that $M^{min}_t$ only depends on that connected component
of the regular values of $\mu$ in which $t$ lies, and as a consequence of this 
we can see that the complex structure on $Q_t$ is the same as on $Q_{t'}$ if
the interval $[t,t']$ does not contain any critical value of $\mu$. 
We have as a conclusion the following theorem:

\begin{theorem} At a regular level $t \in \R$ of the moment map $\rm \mu$, 
we have the K\"ahler quotient $Q_t=\mu^{-1}(t)/S^1$ which is a 
compact K\"ahler orbifold with $M^{min}_t$ as a holomorphic $\C^*$ principal 
 orbibundle on it. Moreover $M^{min}_t$ and the complex structure on $Q_t$
only depend on that connected component of the regular values of $\mu$ where
$t$ lies. 
\label{complex}   
\end{theorem}

It follow from the above theorem that there is a discrete family of complex
orbifolds which arise from the above construction. Moreover, at each level
we get a K\"ahler form on the corresponding complex orbifold. 
The evolution of the
different K\"ahler quotients has been well investigated (e.g. in the papers
\cite{duistermaat-heckman}, \cite{guillemin-sternberg}, cf also
\cite{thaddeus3} and \cite{brion-procesi}). We can summarize these results in the following theorem:

\begin{theorem} The K\"ahler quotients $Q_t$ and $Q_{t'}$ are biholomorphic
if the interval $[t,t']$ does not contain a critical value of the moment
map. They are related by a blowup followed by a blow-down if the interval
$[t,t']$ contains exactly one critical point $c$ different 
from the endpoints. To be more precise, $Q_t$ blown up along the union of
submanifolds 
$\bigcup_{\mu(N_\alpha)=c}P_w(E^-_{\alpha})$ is isomorphic to
$Q_{t'}$ blown up along $\bigcup_{\mu(N_\alpha)=c} P_w(E^+_{\alpha})$ and
in both cases the exceptional divisor is 
$\bigcup_{\mu(N_\alpha)=c} P_w(E^+_{\alpha})\times_{N_\alpha}P_w(E^-_{\alpha})$
the fibre product of weighted projective bundles over $N_\alpha$.

Moreover, in a connected component of the regular values of $\mu$ 
the cohomology classes of the K\"ahler forms $\omega_t(Q_t)$ depend
linearly on $t$ according to the formula:
$$[\omega_t(Q_t)]-[\omega_{t'}(Q_{t'})]=(t-t')c_1(M_t^{min})=
(t-t')c_1(M_{t'}^{min}), $$
 where $c_1$ is the first Chern class of the $\C^*$ principal bundle.
\label{quotients}
\end{theorem}

\subsection{Symplectic cuts}
\label{rizsa3}

Now let us recall the construction of the symplectic cut we need
(see \cite{larman} and also \cite{edidin-graham} for the algebraic case), 
first in the symplectic and second in the K\"ahler
category. 

If $(M,\omega)$
is a symplectic manifold with a Hamiltonian and semi-free
$S^1$ action and proper
moment map $\mu$ with absolute minimum $0$,
then we can define the symplectic cut of $M$ at the regular level $t$ by a
symplectic quotient construction as follows. 

We let $S^1$ act on the
symplectic manifold $M\times \C$ (where the symplectic structure
is the product of the symplectic structure on $M$ and the standard
symplectic structure on $\C$) by acting on the first factor according
to the above $S^1$ action and on the second factor
by the standard multiplication.
This action is clearly Hamiltonian with proper moment map
$\mu+\mu_{\C}$, where
$\mu_\C$ is the standard moment map on $\C$: $\mu_C(z)={\left| z\right|}^2$.

Now if $t$ is a regular value of the moment map $\mu+\mu_\C$, such
that $S^1$ acts with finite stabilizers on $M_t=\mu^{-1}(t)$ (i.e.
$M_{t}/S^1$ gives a symplectic orbifold), then the
symplectic quotient $\bar{M}_{\mu<t}$ defined by
$$\bar{M}_{\mu<t}=\{(m,w)\in M\times \C: \mu(m)+\left| w\right|^2=t\}$$
will be a symplectic compactification of the symplectic manifold
$M_{\mu<t}$ in the sense that $$\bar{M}_{\mu<t}=M_{\mu<t}\cup
Q_{t},$$ and the inherited symplectic structure on $\bar{M}_{\mu<t}$
restricted to $M_{\mu<t}$ coincides with its original symplectic structure.
Moreover, if we restrict this structure onto $Q_t$ it coincides with
its quotient symplectic structure.

Now suppose that we are given a K\"ahler manifold $(M,I,\omega)$ and a
holomorphic $\C^*$ action on it, such that the induced $S^1\subset \C^*$ action
preserves the K\"ahler form and is semi-free and 
Hamiltonian with proper moment map. With 
these extra structures the symplectic cut construction will give us 
$\bar{M}_{\mu<t}$ a compact K\"ahler
orbifold with a $\C^*$ action, such that $\bar{M}_{\mu<t}\setminus Q_t$ is
symplectomorphic to $M_{\mu<t}$ as above and furthermore is biholomorphic 
to $\C^*(M_{\mu<t})$, the union of $\C^*$-orbits intersecting $M_{\mu<t}$.
(This is actually an important point, as it shows
that $M_{\mu<t}$ is {\em not} K\"ahler isomorphic to
 $\bar{M}_{\mu<t}\setminus Q_t$, cf \cite{larman}). 
We can collect all these results into the next theorem:

\begin{theorem} The symplectic cut $\bar{M}_{\mu<t}=M_{\mu<t}\cup Q_t$ 
as a symplectic manifold compactifies the symplectic manifold $M_{\mu<t}$, 
such that the restricted symplectic structure on $Q_t$ coincides with the 
quotient symplectic structure.

Furthermore if $M$ is a K\"ahler manifold with a $\C^*$ action as above, then
$\bar{M}_{\mu<t}$ will be a K\"ahler orbifold with a $\C^*$ action, such
that $Q_t$ with its quotient complex structure is a codimension $1$ 
complex suborbifold of
$\bar{M}_{\mu<t}$ whose complement is equivariantly biholomorphic to $\C^*(M_{\mu<t})$ with
its canonical $\C^*$ action.
\label{cut}
\end{theorem} 

\begin{remark} Note that if $t$ is higher than the highest critical value
(this implies that we have finitely many of them),
then $\C^*(M_{\mu<t})=M$ is the whole space, 
therefore the symplectic
cutting in this case gives a holomorphic compactification of $M$
itself. The compactification is 
$\bar{M}_{\mu < t}$, which is equal to the quotient of 
$(M\times {\bf C} - N \times \{0 \})$ by the action of ${\bf
C}^{\ast}$, where $N$ is the downward Morse flow.
This is the compactification we shall examine here for the case of $\M$, the
moduli space of stable Higgs bundles with fixed determinant of degree $1$.
\end{remark}

\section{The moduli of Higgs bundles $\M$} 
\label{higgs}

\begin{notation} Let
\begin{itemize}
\item
$\Sigma$ be a closed Riemann surface of genus
$g>1$,
\item
$\Lambda$ a fixed line bundle on $\Sigma$ of degree 1,
\item
$\N$ the moduli space of rank $2$ stable bundles with
determinant $\Lambda$,
\item
$\M$
denote the moduli space of stable Higgs pairs
$(E,\Phi)$, where $E$ is a rank $2$ vector bundle on $\Sigma$ with
$\det E= \Lambda$ and $\Phi \in H^0(\Sigma, End_0 E\otimes K_{\Sigma})$.
\end{itemize}
\end{notation}

\begin{remark} For the terms used above we refer the
reader to \cite{hitchin1} and \cite{simpson}.
\end{remark}

After introducing the space $\M$, Hitchin gave its extensive description
in \cite{hitchin1}, \cite{hitchin2}. 
Here we restate some of his results.

\begin{itemize}
\item 
$\M$ is a noncompact,
smooth manifold of dimension $12g-12$ containing $T^*_\N$ as a dense
open set.

\item
Furthermore $\M$ is canonically a Riemannian manifold with a complete
hyperk\"ahler metric. Thus $\M$ has complex structures parameterized
by $S^2$. One of the complex structures, for which $T^*_\N$ is a complex
submanifold, is distinguished, call it $I$. 
We will only be concerned with this
complex structure here. The others (apart from $-I$) are 
biholomorphic to each other and
give $\M$ the structure of a Stein manifold. From these K\"ahler forms one
can build up a holomorphic symplectic form $\omega_h$ on ($\M,I$). 

\item

There is a map, called the Hitchin map $$\chi:\M\rightarrow
H^0(\Sigma,K_{\Sigma}^2)=\C^{3g-3}$$ defined by
$$(E,\Phi)\mapsto \det \Phi.$$ The Hitchin map is proper and an algebraically
completely integrable Hamiltonian system with respect to the holomorphic
symplectic form $\omega_h$, with generic fibre a Prym variety corresponding
to the spectral cover of $\Sigma$ at the image point.

\item

Let $\omega$ denote the K\"ahler form corresponding to the complex structure
$I$. There is a holomorphic $\C^*$ action on $\M$ defined by 
$(E,\Phi)\mapsto (E,z\Phi)$. The restricted action of $S^1$ defined by
$(E,\Phi)\rightarrow (E,e^{i\theta}\Phi)$ is isometric and indeed
Hamiltonian with proper moment map $\mu$. The function 
$\mu$ is a perfect Morse function,
moreover:

$\mu$ has $g$ critical values: an absolute minimum $c_0=0$ and
$c_d=(d-\frac{1}{2})\pi$, where $d=1,...,g-1$.

$\mu^{-1}(c_0)=\mu^{-1}(0)=N_0=\N$
is a non-degenerate critical manifold of index $0$.

$\mu^{-1}(c_d)=N_d$ is a non-degenerate critical manifold of index
$2(g+2d-2)$ and is diffeomorphic to a $2^{2g}$-fold cover of the
$(2g-2d-1)$-fold symmetric product $S^{2g-2d-1}(\Sigma)$.

\item

The fixed point set $S$ of the involution $\sigma(E,\Phi)=(E,-\Phi)$ is
the union
of $g$ complex submanifolds of $\M$ namely, 
$$S=\N\cup \bigcup_{d=1}^{g-1}F_{d},$$
where $F_d$ is the total space of a vector bundle $F_d$ over $Z_d$. Moreover
$F_d$ is a complex submanifold of dimension $3g-3$.

\end{itemize}

Using an algebraic point of view Nitsure in \cite{nitsure} could prove:

\begin{theorem}[Nitsure] $\M$ is a quasi-projective variety.
\end{theorem}

The main aim of this paper is to examine in certain sense the
canonical compactification of $\M$.

\begin{example} Unfortunately, even when $g=2$ the moduli space $\M$ is 
already $6$ dimensional, too big to serve as an instructive example. 
We rather choose $\M_{toy}$, the moduli space of stable 
parabolic Higgs bundles on 
$\Proj^1$, with four marked points,  
 in order to show how
our later constructions work. (These moduli spaces were considered by Yokogawa  \cite{yokogawa}.)  
We choose this example because
it is a complex surface, and can be constructed explicitly. 

We fix four distinct points on $\Proj^1$ and denote by $\Proj^1_4$ the 
corresponding complex orbifold. 
Let $P$ be the elliptic curve corresponding to $\Proj^1_4$.
Let $\sigma_P$ be the involution $\sigma_P(x)=-x$ on $P$. 
Thus, $P/\sigma_P$ is just 
the complex orbifold $\Proj^1_4$. The four fixed
points of the involution $x_1,x_2,x_3,x_4\in P$ correspond to the four marked
points on $\Proj^1_4$. Furthermore, let $\tau$ be the involution 
$\tau(z)=-z$ on $\C$. 

Consider now the quotient space 
$(P\times \C)/(\sigma_P\times \tau)$. 
This is a complex orbifold of dimension $2$
with four isolated $\Z_2$ quotient singularities at the points 
$x_i\times 0$. Blowing up these singularities we get a smooth complex 
surface $\M_{toy}$ with four exceptional divisors $D_1,D_2,D_3$ and $D_4$.  
Moreover the map $\chi:(P\times \C)\rightarrow \C$ sending 
$(z, x)\mapsto z^2$, 
descends to the quotient $(P\times \C)/(\sigma_P\times \tau)$ and sending 
the exceptional divisors to zero one obtains a map 
$\chi_{toy}:\M_{toy}\rightarrow\C$, 
with generic fibre $P$. 
The map $\chi_{toy}$ will serve as our toy Hitchin map.

Moreover there is a $\C^*$ action on $\M_{toy}$, coming from the standard
action on $\C$. 
The fixed point set of $S^1\subset\C^*$ has five
components: one is 
$\N_{toy}\subset \M_{toy}$ (the moduli space of stable parabolic
bundles on $\Proj_4^1$) which is the proper transform of 
$(P\times 0)/(\tau\times \sigma_P)=
P^1_4\subset (P\times \C)/(\sigma_P\times \tau)$ in $\M_{toy}$. 
The other four components consist of single points 
$\tilde{x_i}\in D_i,\ i=1,2,3,4$.

The fixed point set of the involution $\sigma:\M_{toy}\rightarrow \M_{toy}$
has five components, one of which is $\N_{toy}$, the others $F_i$ 
are the proper
transforms of the sets 
$(x_i\times \C)/(\sigma_P\times \tau)
\subset (P\times \C)/(\sigma_P\times \tau).$ 
\end{example}

\section{The nilpotent cone $N$}
\label{nilpotent}

 The results in the previous section show that the K\"ahler manifold 
 $(\M,I,\omega)$ is equipped with a $\C^*$ action which restricts
to an $S^1\subset \C^*$ action which is semi-free and Hamiltonian with proper
moment map $\mu$. Moreover, $0$ is an absolute minimum for $\mu$. Therefore
we are in the situtation described in 
Section~\ref{rizsa}. In the following sections
we will apply the ideas developed there to our situation and deduce
important properties of the spaces $\M$, $Z$ and $\cM$.

We saw in Theorem~\ref{stratification} that the 
downward Morse flow is a deformation retract of $\M$, so it is
responsible for the topology, and as such it is an important object. 
On the other hand we will prove that
the downward Morse flow coincides with the nilpotent cone.

\begin{definition}
The nilpotent cone is the preimage of zero of the Hitchin map
$N=\chi^{-1}(0)$.
\end{definition}

The name `nilpotent cone' was given by Laumon, to emphasize the analogy with 
the nilpotent cone in a Lie algebra. 

In our context this is the most important fibre of the Hitchin map,
and the most singular one at the same time. We will show below that the
nilpotent cone is a central notion in our considerations. 

Laumon in \cite{laumon} investigated the nilpotent cone in a much more
general context and showed its importance in the Geometric Langlands 
Correspondence. Thaddeus in \cite{thaddeus1} concentrated on our case, 
and gave the exact description of the nilpotent cone. In what follows
we will reprove some of their results.  

The following assertion was already stated in \cite{thaddeus1} 
which will turn out to be crucial in some of our considerations.

\begin{theorem} The downward Morse flow coincides with the nilpotent
cone.
\label{morse}
\end{theorem}

\begin{proof} As we saw in Theorem~\ref{stratification} the downward Morse flow can be identified with 
the set of points in $\M$ whose $\C^*$ orbit
is relatively compact in $\M$.

Since the nilpotent cone is invariant under the $\C^*$ action
and compact ($\chi$ is proper) we immediately get that the nilpotent
cone is a subset of the downward Morse flow.

On the other hand if a point in $\M$ is not in the nilpotent cone then the
image of its $\C^*$ orbit by the Hitchin map is a line in $\C^{3g-3}$, 
therefore cannot be relatively compact.
\end{proof}

Laumon's  main
result is the following assertion (cf. Theorem $3.1$ in \cite{laumon}), 
which we prove in our case:

\begin{corollary}[Laumon] The nilpotent cone is a Lagrangian subvariety 
of $\M$
with respect to the holomorphic symplectic form $\omega_h$.
\label{lagrange}
\end{corollary}

\begin{proof} The Hitchin map is a completely integrable Hamiltonian
system, and the nilpotent cone is a fibre of this map, so it is
coisotropic. Therefore it is Lagrangian if and only if its dimension is
$3g-3$.

On the other hand the nilpotent cone is exactly the downward Morse flow and
we can use Hitchin's description of the critical submanifolds in 
\cite{hitchin1}, giving that
the sum of the index and the real dimension of
any critical submanifold is $6g-6$.
We therefore conclude that the complex dimension of the downward Morse flow 
(i.e. the nilpotent cone) is $3g-3$.
\end{proof}

\begin{remark} Nakajima's Proposition $7.1$ in \cite{nakajima} states that
if $X$ is a K\"ahler manifold with a $\C^*$ action and a 
holomorphic symplectic form $\omega_h$ of homogeneity $1$ then the downward
Morse flow of $X$ is Lagrangian with respect to $\omega_h$. Thus Nakajima's 
result and Theorem~\ref{morse} together give an alternative proof of the 
theorem. We prefered the one above for it concentrates on the specific 
properties of $\M$. 
\end{remark}

From the above proof we can see that for higher rank Higgs bundles
Laumon's theorem is equivalent to the assertion that every critical
submanifold contributes to the middle dimensional cohomology, i.e the
sum of the index and the real dimension of
any critical submanifold should always
be half of the real dimension of the corresponding moduli space.

Using the results of \cite{gothen1} one easily shows 
that the above statement also holds 
for the rank $3$ case. 
Gothen could show directly the above statement
for any rank and therefore gave an alternative proof of Laumon's theorem
in these cases \cite{gothen2}.

\begin{corollary} The middle dimensional homology $H_{6g-6}(\M)$ of $\M$ is
freely generated by the homology classes of components
of the nilpotent cone and therefore has dimension $g$.
\label{middle}
\end{corollary}

\begin{proof} We know that each component of $N$ is a projective variety of
dimension $3g-3$. $N$ is a deformation retract of $\M$, therefore
the middle dimensional homology of $\M$ is generated by the homology classes
of the components of $N$. Furthemore, from the Morse picture,
components of $N$ are in a one to
one correspondence with the critical manifolds of $\M$, so 
there are $g$ of them.
The result follows. 
\end{proof}

We finish this section with Thaddeus's description of the nilpotent cone
(see \cite{thaddeus1}, cf. \cite{laumon}).

\begin{theorem} The nilpotent cone is the union of $\N$ and the total
spaces  of vector bundles $E^-_d$ over $N_d$, where $E^-_d$ is the negative
subbundle of $T_\M\mid_{N_d}$.

Moreover, the restricted action of $\C^*$ on $N$ is just the inverse
multiplication 
on the fibres. 
\label{th-nil}
\end{theorem}

\begin{proof} This follows directly from Theorem~\ref{stratification} and
Theorem~\ref{morse}, with noting that by Hitchin's description of the weights
of the circle action on $T_\M\mid N_d$ in the proof of 
Proposition $7.1$ of \cite{hitchin1},
we have that there is only one negative weight. Therefore the $\beta$-fibration
of Theorem~\ref{stratification} is a vector bundle in this case. 
The result follows. 
\end{proof}

\begin{remark} From the description of $E^-_d$ in \cite{thaddeus1} and that of
$F_d$, a component of the fixed point set
of the involution $\sigma(E,\Phi)=(E,-\Phi)$,
in \cite{hitchin1}, one obtains the remarkable fact 
that the vector bundle $E^-_d$ is actually dual to $F_d$. 
\end{remark}

\begin{example} In our toy example we have the elliptic fibration 
$\chi_{toy}:\M_{toy}\rightarrow \C$, with the only singular fibre 
$N_{toy}=\chi_{toy}^{-1}(0)$, the toy nilpotent cone. We have now
the decomposition $$N_{toy}=\N_{toy}\cup \bigcup^4_{i=1}D_i,$$ where
we think of $D_i$ as the closure of $E_i$, the total space of the trivial
line bundle on $\tilde{x_i}$. 
 
The possible singular fibres of elliptic fibrations
have been 
classified by Kodaira (cf. \cite{barth-peters-van}, p. 150). According
to this classification $N_{toy}$ is of type $I_0^*(\tilde{D}_4)$.     
\end{example}

\section{The highest level K\"ahler quotient $Z$}
\label{kaehler}

In this section we apply the ideas of Subsection~\ref{rizsa2} to our 
situation. 

\begin{definition} Define for every non negative $t$ the K\"ahler
quotient $$Q_t=\mu^{-1}(t)/S^1.$$ 
As the complex structure of the K\"ahler quotient
depends only on the connected component of the regular values of $\mu$, we
can define $Z_d=Q_t $ for $c_{d}<t<c_{d+1}$ as a complex orbifold (we take
$c_g=\infty$). Similarly, we define $X_{Z_d}$ to be $\M^{min}_t$ for 
$c_{d}<t<c_{d+1}$.

For simplicity let the highest level quotient $Z_{g-1}$
be denoted by $Z$ and the corresponding $\C^*$ principal bundle $X_{Z_{g-1}}$ 
by $X_{Z}$. 
\end{definition}

In the spirit of Theorem~\ref{quotients} we have the following

\begin{theorem}$Z_d$ is a complex orbifold with only
$\Z_2$-singularities, where the singular
locus is diffeomorphic to some union of projectivized vector bundles
$P(F_i)$:
$$\sing(Z_d)=\bigcup_{0<i\leq d} P(F_i), $$
where $F_i\subset\M$ is the total space of a vector bundle over $N_i$ and 
is a component of the fixed point set of the involution 
$\sigma(E,\Phi)=(E,-\Phi)$.  

\label{Zorbi}
\end{theorem}

\begin{proof} 
The induced action of $S^1$ on $\C^{3g-3}$ by the Hitchin map
is multiplication by $e^{2i\theta}$ so an orbit of $S^1$ on
$\M\setminus N$ is a non trivial double cover 
of the image orbit on $\C^{3g-3}$.
On the other hand by Thaddeus' description of $N$ (Theorem~\ref{th-nil}) 
it is clear 
that if a point of $N$
is not a fixed point of the circle action, then the stabilizer is trivial at
that point.

Summarizing these two observations we obtain that if a point of $\M$ is
not fixed by $S^1$, then its stabilizer is either trivial or $\Z_2$ . The
latter case occurs exactly at the fixed point set of the involution $\sigma$.
The statement now follows from Theorem~\ref{quotients}.
\end{proof}

\begin{proposition}
$Z_d$ and $Z_{d+1}$ are related by a blowup following by a blowdown.

Namely, $Z_d$ blown up along $P(E^-_d)$ is the
same as the singular quotient $Q_{c_d}$ blown up along $N_d$ (its singular
locus),
which in turn gives $Z_{d+1}$ blown up at $P(E^+_d)$.

Moreover, this birational equivalence is an isomorphism outside an analytic set
of codimension at least $3$.
\label{birational}
\end{proposition}

\begin{proof} The first bit is just the 
restatement of Theorem~\ref{quotients} in
our setting.

The second part follows because $\dim(P(E_d^-))=3g-3-1<6g-6-2$ and
$\dim(P(E^+_d))=3g-3+2g-2d-1-1<6g-6-2$ for $g>1$.
\end{proof}

\begin{corollary} $Z=Z_{g-1}$  is birationally equivalent to $P(T^*_\N)=Z_0$.
Moreover this gives an isomorphism in codimension $>2$.
\label{Zjeno}
\end{corollary}

\begin{proof} Obviously $X_{Z_0}$ is $T^*_\N$, and therefore by 
Theoreom~\ref{complex} $Z_0$ is isomorphic 
to the projectivized cotangent bundle $P(T^*_\N)$. The
statement follows from the previous theorem. 
\end{proof}

\begin{corollary} $Z$ has  Poincar\'e polynomial
$$ P_t(Z)=\frac{t^{6g-6}-1}{t^2-1}P_t(\N) +
\sum _{i=1}^{g-1}\frac{t^{6g-6}-t^{2g-4+4i}}{t^2-1}P_t(N_i),$$
where $N_i$ is a $2^{2g}$-fold cover of $S^{2i-1}\Sigma$.
\label{PZ}
\end{corollary}

\begin{proof} One way to derive this formula is through Kirwan's formula
in \cite{kirwan}. We use the above blowup, blowdown picture instead. 
This approach is due to Thaddeus, see \cite{thaddeus2}.

Applying the formula in \cite{griffiths-harris},p.605 twice we get that
$$P_t(Z_{d+1})-P_t(Z_d)=P_t(PE^+_d)-P_t(PE^-_d).$$

On the other hand for a projective bundle on a manifold
$P\rightarrow M$ with fiber $\Proj^n$ one has
(cf. \cite{griffiths-harris} p.606)
$$P_t(P)=\frac{t^{2n+2}-1}{t^2-1}P_t(M).$$

Hence the formula follows.
\end{proof}

\begin{remark} All the Poincar\'e polynomials on the right hand side
of the above formula have been calculated.
For $P_t(\N)$ see e.g. \cite{atiyah-bott} for $P_t({N_d})$ see
\cite{hitchin1}.
\end{remark}

We will determine the Picard group of $Z$ exactly. First we
define some line bundles on several spaces.

\begin{notation} Let 

\begin{itemize} 
\item $\calL_\N$ denote the ample generator 
of the Picard group of $\N$
(cf. \cite{drezet-narasimhan}) 

\item
$\calL_{PT^*_\N}$ be its pullback 
to $PT^*_\N$, 
\item
$\calL_Z$ 
denote the corresponding line bundle on $Z$ (cf. Corollary~\ref{Zjeno}).

\item 
$L_{PT^*_\N}$ be the dual of the tautological line bundle on $PT^*_\N$,

\item 
$L_Z=X^*_Z\times_{\C^*}\C$ denote the corresponding line orbibundle
on $Z$. 
\end{itemize}
\end{notation}

\begin{corollary} $\Pic(Z)$, the Picard group of $Z$, is of rank $2$ over
$\Z$ and is freely generated by $\calL_Z$ and $L_Z$.
\label{picard}
\end{corollary}

\begin{remark} The Picard group of $Z$ is the group of invertible sheaves on
$Z$. As the singular locus of $Z$ has codimension $\geq 2$, this group can
be thought of as the group of holomorphic line orbibundles on $Z$. Namely, in
this case the restriction of a holomorphic line orbibundle to 
$Z\setminus \sing(Z)$ gives a one-to-one correspondence between holomorphic
line orbibundles on $Z$ and holomorphic line bundles on $Z\setminus \sing(Z)$,
by the approriate version of Hartog's theorem.
\end{remark}

\begin{proof} It is well known that $\Pic(\N)$ is freely generated by one ample
line bundle $\calL_\N$ therefore is of rank $1$ (cf. \cite{drezet-narasimhan}).
Thus $\Pic(P(T^*_\N))$
is  of rank $2$ and freely generated by $\calL_{PT^*_\N}$
the pullback of $\calL_\N$ and the dual of the 
tautological line bundle $L_{PT^*_\N}$.
From Corollary~\ref{Zjeno} $\Pic(Z)$ is
isomorphic with $\Pic(P(T^*_\N))$ therefore is of rank $2$, and freely
generated by $\calL_Z$ and $L_Z$, where $\calL_Z$ is isomorphic
to $\calL_{PT^*_\N}$ and $L_Z$ is isomorphic to $L_{PT^*_\N}$ outside
the codimension $2$ subset of Corollary~\ref{Zjeno}.
\end{proof}

\begin{definition} A contact structure on a compact complex orbifold
$Z$ of complex dimension $2n-1$ is given by the following data: 

\begin{enumerate}
\item
a contact line orbibundle $L_Z$ such that $L_Z^n=K_Z^{-1}$, where  $K_Z$ is
the line orbibundle of the canonical divisor of $Z$,

\item

a complex contact form $\theta\in H^0(Z,\Omega^1(Z)\otimes L_Z)$ a 
holomorphic $L_Z$ valued $1$-form, such that 
\begin{eqnarray}
0\neq \theta\wedge (d\theta)^{n-1}\in H^0(Z,\Omega^{2n-1}(Z)\otimes K_Z^{-1})=
H^0(Z,\calO_Z)=\C 
\label{degeneralt}
\end{eqnarray}
is a nonzero constant.
\end{enumerate} 
\end{definition}
 
\begin{theorem} There is a canonical holomorphic contact structure
on $Z$ with contact line orbibundle $L_Z$.
\label{contact}
\end{theorem}

\begin{proof} This contact structure can be created by the construction of 
Lebrun as in \cite{lebrun} Remark $2.2$. 
We only have to note that the holomorphic 
symplectic form $\omega_h$ on $\M$ is of homogeneity $1$. 

The construction goes as follows. If $\pi:X^*_Z\rightarrow Z$ denotes the
canonical projection of the $\C^*$ principal orbibundle $X^*_Z$ the
dual of $X_Z$, then $\pi^*(L_Z)$ is canonically trivial with the 
canonical section having homogeneity $1$. Thus in order to give a complex
contact form $\theta\in H^0(Z,\Omega^1(Z)\otimes L_Z)$ it is sufficient
to give a $1$-form $\pi^*\theta$ on $X^*$ of homogeneity $1$. This can
be defined by $\pi^*\theta=i(\xi)\omega_h$, where $\xi\in H^0(\M,T_\M)$
is the holomorphic vector field generated by the $\C^*$ action. The 
non-degeneracy condition (\ref{degeneralt}) is exactly equivalent to requiring
that the closed holomorphic $2$ form $\omega_h$ satisfy $\omega_h^n\neq 0$. 
This is the case as $\omega_h$ is a holomorphic symplectic form. 

The result follows.
\end{proof}

We will be able to determine the line orbibundle $L_Z$ explicitly. 
For this, consider the Hitchin map $\chi:\M\rightarrow \C^{3g-3}$. 
As it is equivariant
with respect to the $\C^*$ action, $\chi$ induces a map 
$$\bar{\chi}:Z\rightarrow \Proj^{3g-4}$$
on Z. The generic fibre of this map is easily seen to be the Kummer
variety corresponding to the Prym variety (the Kummer variety of an
Abelian variety is the quotient of the Abelian variety by the involution
$x\rightarrow -x$), the generic fibre of the Hitchin map. 
Thus we have proved 

\begin{lemma} There exists a map $\chi:Z\rightarrow \Proj^{3g-4}$ the reduction
of the Hitchin map onto $Z$, for which the generic fibre is a Kummer variety.
\label{Zchi}
\end{lemma}

\begin{remark} This observation was already implicit in Oxbury's thesis 
(cf. $2.17$a of \cite{oxbury}). \end{remark}

The following theorem determines the line bundle $L_Z$ in terms of the
Hitchin map.

\begin{theorem} $L^2_Z=\bar{\chi}^*{\cal H}_{3g-4}$ where 
${\cal H}_{3g-4}$ is the hyperplane bundle on $\Proj^{3g-4}$.
\label{feco}
\end{theorem}

\begin{proof} We understand from Corollary~\ref{picard} that
$\bar{\chi}^* {\cal H}_{3g-4}={\calL}_Z^k\otimes L_Z^l$ for some 
integers $k$ and $l$.

We show that $k=0$. For this consider the pullback of $\calL_Z$ onto
$\M\setminus N$ the total space of the $\C^*$ principal orbibundle $X^*_Z$. 
This line orbibundle extends to $\M$ as $\calL_\M$
and restricts to $T^*_\N$ as the pullback
of $\calL_{PT^*_\N}$ by construction. $c_1(\calL_\M)$ is not trivial when
restricted to $\N$ 
(namely it is $c_1(\calL_\N)$, since this bundle is ample) 
therefore is not trivial when
resticted to a generic fibre of the Hitchin map. We can deduce
that $c_1(\calL_Z)$ is not trivial on the generic fibre of $\bar{\chi}$.

However $L_Z$ restricted to a generic fiber of
$\bar{\chi}$ can be described as follows. Let this Kummer variety
be denoted by $K$, the corresponding Prym variety by $P$. Form
the space $P\times \C^*$, the trivial $\C^*$ principal bundle
on $P$ and quotient it out by the involution
$\tau(p,z)=(-p,-z)$. The resulting space is easily seen to be
the $\C^*$ orbit of the Prym $P$ in $\M$, therefore the total space
of the $\C^*$ principal orbibundle $L^*_Z\setminus (L^*_Z)_0$ on $K$. 
Hence $L_Z^2$ is the trivial line orbibundle on $K$. Thus $c_1(L_Z\mid_K)=0$.

Now
$\bar{\chi}^*{\cal H}_{3g-4}$ is trivial on the Kummer variety.
Hence the assertion
$k=0$.

The rest of the proof will follow the lines of Hitchin's proof of Theorem
$6.2$ in \cite{hitchin2}. We show that $l=2$.

The sections of $L_Z$ can be 
identified with holomorphic functions
homogeneous of degree $2$ on the $\C^*$ principal
orbibundle $X_Z=L_Z^*\setminus (L_Z^*)_0$. As $N$ is of codimension
$\geq 2$ such functions extend to $\M$. Since the Hitchin map is proper,
these functions are constant on the fibers of the Hitchin map, therefore
are the pullbacks of holomorphic functions on $\C^{3g-3}$ of
homogeneity $1$ which can be identified with the holomorphic sections of
the hyperplane bundle ${\cal H}_{3g-4}$ on $P(\C^{3g-3})=\Proj^{3g-4}$.
\end{proof}

\begin{corollary} If $n$ is odd, there are natural isomorphisms

$$H^0(Z,L_Z^n)\cong H^0(\N,S^n{T_\N})\cong0,$$ whereas if $n$ is even, then
$$ H^0(Z,L_Z^n)\cong H^0(\N,S^n{T_\N})
\cong H^0(\Proj^{3g-4},{\cal H}^{\frac{n}{2}}_{3g-4}).$$
\label{H0Z}
\end{corollary}

\begin{proof} We show that $H^0(Z,L_Z)\cong H^0(\N,S^n(T_\N))$ for every
$n$, the rest of the theorem will follow from Theorem $6.2$ of 
\cite{hitchin2}. 

By Proposotion~\ref{birational} we get that 
$H^0(Z,L^n_Z)\cong H^0(PT^*_\N,L^n_{PT^*_\N})$. Let $\pi:PT^*_\N\rightarrow\N$
denote the projection. It is well known that the Leray spectral sequence
for $\pi$ degenerates at the $E^2$ term. Moreover, we have that 
$R^i\pi_*(L^n_{PT^*_\N})=0$ if $0<i<{3g-4}$ 
(cf. \cite{hartshorne} Theorem $5.1$b). Therefore 
$H^0(PT^*_\N,L^n_{PT^*_\N})\cong H^0(\N,\pi_*(L^n_{PT^*_\N}))$. Finally 
the sheaf $\pi_*(L^n_{PT^*_\N})$ is $S^n(T_\N)$, which proves the statement.
\end{proof}

We can moreover determine the first cohomology group corresponding to
the infintesimal deformations of the holomorphic contact structure on $Z$
and can interpret it in a nice way.

\begin{corollary} There are canonical isomorphisms 
$$H^1(Z,L_Z) \cong (H^1(\M,{\calO}_{\M}))_1\cong H^1(\N,T_\N)
\cong H^1(\Sigma, K_\Sigma^{-1}),$$
where $(H^1(\M,{\calO}_{\M}))_1\subset H^1(\M,{\calO}_{\M})$ is 
the vector space of elements of $H^1(\M,{\calO}_{\M})$ homogeneous of
degree $1$.
\label{H1Z}
\end{corollary}

\begin{proof} We may use the cohomological version of Hartog's theorem 
(cf. \cite{scheja}) to show that $H^1(Z,L_Z)\cong H^1(PT^*_\N,L_{PT^*_\N})$,
as $Z$ and  $PT^*_\N$ are isomorphic on an analytic set of codimension $\geq 3$
(cf. Proposition~\ref{birational}).

The proof of the other isomorphisms can be found in \cite{hitchin3}.
\end{proof}

\begin{remark} We can interpret this result as saying that the deformation of
the complex structure on $\Sigma$ corresponds to the deformation of complex
structure on $\N$, to the deformation of holomorphic contact structure on $Z$
(cf. \cite{lebrun}) and to the deformation of the holomorphic symplectic
structure of homogeneity $1$ on $\M$.
\end{remark}

As an easy corollary of the above we note the following

\begin{corollary} The line orbibundle $L_Z$ is nef but neither trivial
nor ample.
\label{triv}
\end{corollary}

\begin{proof} $L_Z$ is certainly not ample since $c_1(L_Z)$ is trivial on the
Kummer variety.

On the other hand $L^2_Z$ being the pullback of an ample bundle is not
trivial and is nef itself,
hence the result. 
\end{proof}

The next theorem will describe the inherited K\"ahler structures of $Z$.
Considering the one-parameter family of
K\"ahler quotients $Q_t$, $t>c_{g-1}$ we
get a one-parameter family of K\"ahler forms $\omega_t$ on $Z$.
Theorem $1.1$ from \cite{duistermaat-heckman} gives the
following result for our case (cf. Theorem~\ref{quotients}). 

\begin{theorem}[Duistermaat,Heckman]
The complex orbifold $Z$ 
has a one-parame\-ter family of K\"ah\-ler forms $\omega_t$,
$t>c_{g-1}$ such that 
$$[\omega_{t_1}(Z)]-[\omega_{t_2}(Z)]=(t_1-t_2)c_1(L_Z)$$
where $t_1,t_2>c_{g-1}$ and $[\omega_t]\in H^2(Z,\R)$ is the cohomology
class of $\omega_t$.
\label{csalad}
\end{theorem}

Many of the above results will help us to prove the following theorem.

\begin{theorem} $Z$ is a projective algebraic variety.
\label{proj}
\end{theorem}

\begin{proof} 
By the Kodaira embedding theorem for orbifolds (cf. \cite{bailey}) we have only
to show that $Z$ with a suitable K\"ahler form is a Hodge orbifold, i.e.
the K\"ahler form is integer. For this to see we show that the
K\"ahler cone of $Z$ contains a subcone, which is open in $H^2(Z,\R)$.
This is sufficient since such an open subcone should contain an
integer K\"ahler
form i.e. a Hodge form.

Since Corollary~\ref{picard} shows that $\Pic_0(Z)$ is trivial, 
by Corollary~\ref{triv} we see that $c_1(L_Z)\neq 0$. Therefore the 
 previous theorem exhibited a half line in the K\"ahler cone of $Z$. Thus
to find an open subcone in the $2$ dimensional vector space $H^2(Z,\R)$
(Corollary~\ref{picard}) it is sufficient to show that this line 
does not go through
the origin or in other words $c_1(L)$ is not on the line.
But this follows from Corollary~\ref{triv}, because $L$ being not ample
$c_1(L)$ cannot contain a K\"ahler form. Hence the result.
\end{proof}

\begin{remark} We see from this proof that $c_1(L_Z)$ lies on the closure
of the K\"ahler cone, thus $L_Z$ is nef. This reproves a statement of
Corollary~\ref{triv}.
\end{remark}

\begin{example}  In the case of the toy example the lowest level
K\"ahler quotient $Z_0$ is    
the projectivized cotangent bundle 
$PT^*_{\N_{toy}}$ of $\N_{toy}$, which is isomorphic to 
$\N_{toy}=\Proj^1$, and
the blowups and blowdowns 
add the four marked points to $\Proj^1$.  
Therefore $Z_{toy}$ is isomorphic to the orbifold $\Proj^1_4$, where
the marked points correspond to the fixed point set of the involution 
$\sigma$, namely these are the projectivized bundles $PF_i$, i.e. points.

Moreover the $\C^*$ principal orbibundle $X_{Z_{toy}}$ on $\Proj^1_4$ has the
form 
$$X_{Z_{toy}}=(P\times \C^*)/(\sigma_P\times \tau).$$ 

Thus in the toy example, not
like in the ordinary Higgs case,  
we have $c_1(L_{Z_{toy}})=0$. This latter assertion can be seen using 
\ref{feco} and noting that the target of the reduced toy Hitchin map
$\bar{\chi}_{toy}:Z_{toy}\rightarrow \Proj^0$ is a point. 

There is an other difference, namely the Picard group of $Z_{toy}$ is of 
rank $1$, because $L_{Z_{toy}}^2$ is the trivial bundle on $Z_{toy}$. 
\end{example}

In the next section we show how to compactify $\M$ by sewing in $Z$ at
infinity.

\section{The compactification $\cM$}
\label{kompakt}

In this section we compactify $\M$ by adding to each non-relatively 
compact
$\C^*$ orbit an extra point i.e. sewing in $Z$ at infinity. 
Another way of saying
the same is to glue together $\M$ and $E$ the total space of $L_Z$ along the
$\C^*$ principal orbibundle $X^*_Z=E\setminus E_0=\M\setminus N$. 
To be more
precise we use the construction of Lerman, called the symplectic cut
(cf. Subsection~\ref{rizsa3} and \cite{larman}).

Since the complex structure on the K\"ahler quotients depends only on the
connected component of the level, we can make the following definition.

\begin{definition} Let $\bar{\M}_d$ denote
the compact complex orbifold corresponding
to the K\"ahler quotients of $\M\times \C$ by  the product $S^1$ action
$$\bar{M}_{\mu<t}=(\mu+\mu_\C)^{-1}(t)/{S^1},$$ with $c_d<t<c_{d+1}$.

Let $X_{\cM_d}$ denote the corresponding $\C^*$ principal bundle on $\cM_d$.
For simplicity we let $\bar{\M}$ denote $\bar{\M}_{g-1}$ and $X_\cM$ denote
$X_{\cM_{g-1}}$.
\end{definition}

As a consequence of the construction of symplectic cutting we have the
following theorem (cf. Theorem~\ref{cut})

\begin{theorem} The compact orbifold
$\cM=\M\cup Z$ is a compactification of $\M$ such
that $\M$ is an open complex submanifold and $Z$ is a codimension one 
suborbifold, i.e. a divisor.

Moreover $\C^*$ acts on $\cM$ extending the action on $\M$ with the points
of $Z$ being fixed.
\label{compactification}
\end{theorem}

In addition to the above we see that we have another
decomposition $\cM=N\cup E$ of $\cM$ into the
nilpotent cone and the total space $E$
of the contact line bundle $L_Z$ on $Z$. Thus the
compactification by symplectic cutting produced the same orbifold
as the two constructions we started this section with.

We start to list the properties of $\cM$. We will mention properties
analogous to properties of $Z$ (these correspond to the fact that both
spaces were constructed by a K\"ahler quotient procedure)
and we will clarify the
relation between $Z$ and $\cM$.

Theorem~\ref{cut} and Theorem~\ref{quotients} give the following result
in our case. 

\begin{theorem} $\cM_d$ is  a compact orbifold. It has a decomposition
$\cM_d=\M_d \cup Z_d$ into an open complex suborbifold $\M_d$ (which
is actually a complex manifold) and a codimension one suborbifold $Z_d$,
i.e. a divisor.
The singular locus of $\cM_d$ coincides with that of $Z_d$:

$$\sing(\cM_d)=\sing(Z_d)=\bigcup_{0<i\leq d} P(F_i)$$
where $F_i$ is a component of the fixed points set of the involution
$\sigma(E,\Phi)=(E,-\Phi)$.

Furthermore, the $\C^*$ action on $\M_d$ extends onto $\cM_d$ with an extra
component $Z_d$ of the fixed point set.
\label{Morbi}
\end{theorem}

We have the corresponding statement of Theorem~\ref{Zjeno}.

\begin{theorem} $\cM=\cM_{g-1}$ is birationally isomorphic to
$\cM_0=P(T^*_\N\oplus{\calO}_\N)$. Moreover, they are isomorphic outside
an analytic subset of codimension at least $3$.
\label{Mjeno}
\end{theorem}

\begin{proof} In a similar manner to the proof 
of Corollary~\ref{Zjeno} we can argue
by noting that $X_{\cM_0}$ is obviously isomorphic to $T^*_\N\oplus 
\calO_\N$ with
the standard action of $\C^*$. Hence indeed 
$\M_0=P(T^*_\N\oplus{\calO}_\N)$. 

By Theorem~\ref{quotients} it is clear that $\cM$ and $\cM_0$ are related
by a sequence of blowups and blowdowns. 
The codimensions of the submanifolds
we apply the blowups are at least $3$ by a calculation analogous to the one
in the proof of Proposition~\ref{birational}.  
\end{proof}

\begin{notation} Let 

\begin{itemize} \item
$\calL_{P(T^*_\N \oplus\calO_\N)}$ denote
the pull back of 
$\calL_\N$ to $P(T^*_\N\oplus \calO_\N)$, 
\item 
$\calL_\cM$ be the corresponding 
line bundle on $\cM$.

\item
$L_{P(T^*_\N\oplus\calO_\N)}$ be the dual of the tautological line bundle
on the projective bundle $P(T^*_\N\oplus \calO_\N)$, 
\item  
$L_\cM=X_\cM\times_{\C^*}\C$ be the corresponding line orbibundle on $\cM$.  
\end{itemize}

\end{notation}

\begin{corollary} $\Pic{\cM}$ is isomorphic to
$\Pic(P(T^*_\N\oplus{\calO}_\N))$ and therefore is of rank $2$ and freely
generated by $L_\cM$ and
$\calL_\cM$.
\label{picM}
\end{corollary}

\begin{proof} 
The previous theorem shows that $\cM$ and $P(T^*_\N\oplus \calO_\N)$ are 
isomorphic outside an analytic subset of codimension at least $2$, thus their
Picard groups are naturally isomorphic. 

However, $\Pic(P(T^*_\N\oplus \calO_\N))$ is freely generated by 
$L_{P(T^*_\N\oplus\calO_\N)}$ and 
$\calL_{P(T^*_\N\oplus\calO_\N)}$. The result follows. 
\end{proof}

\begin{corollary}
The canonical line orbibundle $K_\cM$ of $\cM$ coincides with 
$L_\cM^{-(3g-2)}$. Moreover, $L_\cM$ is the line bundle of the divisor $Z$,
therefore (3g-2)Z is the anticanonical divisor of $\cM$. Finally, $L_\cM$
restricts to $L_Z$ to $Z$. 
\label{vonat}
\end{corollary}

\begin{proof} 
$L_\cM$ by its construction clearly restricts to $L_Z$ on $Z$ and
it is the line bundle of $Z$, as the corresponding statement
 is obviously true for
$P(T^*_\N\oplus \calO_\N).$

The restriction of $K_\cM$ to $\M$ has a non-zero section, namely the 
holomorphic Liouville form $\omega_h^{3g-3}$, thus trivial. Hence
$K_\cM=L_\cM^k$ for some $k\in \Z$.

By the second adjunction formula $K_Z=(K_\cM\otimes [Z])\mid_Z$. The right
hand side equals $L_Z^{-(3g-3)}$ as $L_Z$ is a contact line bundle 
(cf. \ref{contact}). The left hand side can be written as 
$(L_\cM^k\otimes L_\cM)\mid_Z=L_Z^{k+1}$, therefore $k=-(3g-2)$.
\end{proof} 

\begin{lemma} $\chi$ has an extension to $\cM$,  
$$\bar{\chi}:\cM\rightarrow \Proj^{3g-3}$$ such that $\bar{\chi}$ restricted
to $Z$ gives the map of Lemma~\ref{Zchi}.
\label{Mchi}
\end{lemma}

\begin{proof}
We let $\C^*$ act on $\C^{3g-3}\times \C$ by $\lambda(x,z)=
(\lambda^2 x,\lambda z)$. With respect to this action the map 
$(\chi,id_\C):\M\times \C\rightarrow \C^{3g-3}\times \C$ is equivariant. 
Therefore making the symplectic cut 
it reduces to a map $\bar{\chi}:\cM\rightarrow \Proj^{3g-3}$ since the
quotient space 
$(\C^{3g-3}\setminus 0)\times\C /\C^*$ is isomorphic to $\Proj^{3g-3}$. 

The result follows. 
\end{proof}

\begin{remark} In the higher rank case where $\C^*$ acts on the
target space of the Hitchin map with different weights the target space
of the compactified Hitchin map is a weighted projective space.
\end{remark}

\begin{corollary} $L^2_\cM=\bar{\chi}^*{\cal H}_{3g-3}$.
\label{hurka}
\end{corollary}

\begin{proof} Obviously, $\bar{\chi}^*{\cal H}_{3g-3}\mid_\M$ is trivial,
therefore $\bar{\chi}^*{\cal H}_{3g-3}$ is some power of $L_\cM$. 
By \ref{feco} this power is $2$. 
\end{proof}

\begin{theorem}[Duistermaat, Heckman] $\cM$ has a one-parameter family
of K\"ahler forms $\omega_t(\cM)$, $t>c_{g-1}$ such that

$$[\omega_{t_1}(\cM)]-[\omega_{t_2}(\cM)]=(t_1-t_2)c_1(L_{\cM}).$$

Furthermore this one-parameter family of K\"ahler forms restricts to Z as
the one-parameter family of K\"ahler forms of Theorem~\ref{csalad}.
\label{DuHeM}
\end{theorem}

\begin{proof} This is just the application of Theorem~\ref{quotients} and
Theorem~\ref{cut} to our situation.
\end{proof}

\begin{corollary} $\cM$ is a projective algebraic variety.
\label{Mproj}
\end{corollary}

\begin{proof} The argument is the same as for Theorem~\ref{proj}, noting that
by Corollary~\ref{picM} $H^2(\cM,\R)$ is two dimensional and $L_\cM$ is neither
trivial nor ample since $L_\cM\mid Z=L_Z$ (by Corollary~\ref{vonat}) is neither
trivial nor ample (by Corollary~\ref{triv}).
\end{proof}

\begin{remark} 1. The above proof yields that the cohomology class 
$c_1(L_\cM)$ sits
in the closure of the K\"ahler cone of $\cM$, hence $L_\cM$ is nef.

2.  From the previous remark and Corollary~\ref{hurka} 
we can deduce that there is a complete hyperk\"ahler (hence Ricci flat) 
metric on $\M=\cM\setminus Z$, the complement of a nef anticanonical divisor 
of a compact orbifold. 

Therefore our compactification of 
$\M$ is compatible with Yau's problem, which addresses the question: which
non-compact complex manifolds possess a complete Ricci flat metric? 
Tian and Yau in \cite{yau-tian} 
could show that this is the case for the complement of
an ample  anticanonical divisor in a compact complex manifold. (Such manifolds
are called Fano manifolds.) 

The similar statement with ample replaced by nef is an 
unsolved problem.
\end{remark}

\begin{theorem} $\cM$ has Poincar\'e polynomial
$$P_t(\cM)=P_t(\M)+t^2P_t(Z).$$
\label{PM}
\end{theorem}

\begin{proof} We have three different ways of calculating the Poincar\'e
polynomial of $\cM$. The first is through Kirwan's formula in \cite{kirwan},
the second is due to Thaddeus in \cite{thaddeus3}, 
which we used to calculate the Poincar\'e
polynomial of $Z$.

For $\cM$ there is a third method, namely direct Morse theory. All we have to
note is that the $S^1$ action $\cM$ is Hamiltonian with respect to
any K\"ahler form of Theorem~\ref{DuHeM}, and the critical 
submanifolds and corresponding
indices are the same as for $\M$ with one extra critical submanifold $Z$
of index $2$. Hence the result.
\end{proof}

\begin{example} We can describe $\cM_{toy}=\M_{toy}\cup Z_{toy}$ as follows.
As we saw above $\M_{toy}\setminus N_{toy}=X_{Z_{toy}}$. Thus gluing together
$\M_{toy}$ and $E_{toy}$, the total space of the line orbibundle 
$L_{Z_{toy}}$, along $X_{Z_{toy}}$ yields 
$$\cM_{toy}=\M_{toy}\cup_{X_{Z_{toy}}} E_{toy}.$$

One can construct $\cM_{toy}$ directly, as follows. Take 
$\Proj^1=\C\cup \infty$ extending  the involution $\tau$ from $\C$ to 
$\Proj^1$. 
Consider
the quotient $(P\times \Proj^1)/(\sigma_P\times \tau)$. This is a compact
orbifold with eight $\Z_2$-quotient singularities. Blow up four of them
corresponding to $0\in \C$. The resulting space will be isomorphic
to $\cM_{toy}$. The remained four isolated $\Z_2$ quotient singularities will
just be the four marked points of $Z_{toy}\subset \cM_{toy}$, the singular
locus of $\cM_{toy}$. 
\end{example}

We finish this section with a result which gives an interesting relation
between the intersections of the component of the nilpotent cone $N$
in $\M$ 
(equivalently the intersection form on the middle compact cohomology 
$H_{cpt}^{6g-6}(\M)$, cf. Corollary~\ref{middle})  
and the contact structure of $Z$.

\begin{theorem}
There is a canonical isomorphism between the cokernel of $j_\M$ and the
cokernel of $L$, where
$$j_\M:H_{cpt}^{6g-6}(\M)\rightarrow H^{6g-6}(\M)$$ is the
canonical map and
$$L:H^{6g-8}(Z)\rightarrow H^{6g-6}(Z)$$
is multiplication with $c_1(L_Z)$.
\end{theorem}

\begin{proof} We will read off the statement from the following diagram. 

$$
\begin{array}{ccccccccc}
 &  &     &       &       0           &        &         &        & \\
 &  &     &       &      \downarrow   &        &         &        & \\
   &    &     &   & H^{6g-8}(Z) &        &         &        & \\ 
   &    & & &\downarrow     &{\searrow}^L &         &  & \\

0 & \rightarrow & H^{6g-6}_{cpt}(\M)&\rightarrow & H^{6g-6}(\cM)
&\rightarrow &H^{6g-6}(Z)&\rightarrow & 0 \\
 & &  &{\searrow}^j&\downarrow&   &  &  &\\
& & & & H^{6g-6}(\M)&&&&\\  
&&&&\downarrow&&&& \\
&&&&0&&&&\\
\end{array}
$$

We show that both the vertical and horizontal sequences are exact and the
two triangles commute. 

From  the Bialynicki-Birula decomposition of $\cM$ 
we get the short exact sequence of
middle dimensional cohomology groups (recall that $E\subset\cM$ denotes
the total space of the contact line bundle $L_Z$ on Z):
$$0\rightarrow H^{6g-6}_{cpt}(E)\rightarrow H^{6g-6}(\cM)\rightarrow
H^{6g-6}(\M)\rightarrow 0.$$

Applying the Thom isomorphism (which also exists in the orbifold category) 
we can identify $H^{6g-6}_{cpt}(E)$ with $H^{6g-8}(Z)$, this gives
the vertical short exact sequence of the diagram. The horizontal one
is just its dual short exact sequence. 

Finally, the left triangle clearly commutes as all the maps are natural, 
while the right triangle commutes because the original triangle commuted as
above and the canonical map $j_E:H^{6g-6}_{cpt}(E)\rightarrow H^{6g-6}(E)$ 
transforms to $L:H^{6g-8}(Z)\rightarrow H^{6g-6}(Z)$ by
the Thom isomorphism. 

Now the theorem is the consequence of the Butterfly lemma 
(cf. \cite{lang} IV.$4$ p.102), or can be proved
by an easy diagram chasing.

Hence the result follows. 
\end{proof}

\begin{remark} 

1. If the line bundle $L_Z$ was ample then the map $L$ would just be the
Lefschetz isomorphism, and therefore the cokernel would be trivial. In
our case we have $L_Z$ being only nef and the map is not an isomorphism,
the cohomology class of the Kummer variety lying in the kernel. Therefore
the cokernel measures how far is $L_Z$ from being ample. 

2. The cokernel of $j_\M$ measures the degeneracy of the intersection
form on the compactly supported middle dimensional cohomology of $\M$. In this
case also the cokernel is not trivial as the compactly supported
cohomology class of the Prym variety lies in the kernel. This can be seen
by thinking of the Hitchin map as a section of the trivial
rank $3g-3$ vector bundle on $\M$ and considering the ordinary cohomology class
of the Prym variety as the Euler class of this trivial vector bundle, and as
such, the ordinary cohomology class of the Prym variety is trivial indeed.    

3. Notice that the proof did not use any particular property of 
$\M$ therefore the statement is true in the general setting of 
Section~\ref{rizsa}.
\end{remark}

\begin{example} 
1. We determine the dimension of the 
cokernels of the above theorem in the case when
$g=2$, by showing that the intersection form on the compactly supported
middle dimensional cohomology is $0$, i.e. the map $j_\M$ is zero. In the 
previous remark we saw that the compactly supported cohomology class of 
the Prym variety $P$ is in the kernel of $j_\M$. 
It follows from \cite{thaddeus1} that
the compactly supported cohomology of $\N$ and that of $P$ generates the
$2$-dimensional compactly supported middle cohomology 
of $\M$ (cf. Theorem~\ref{middle}). 

On the other hand the Euler
characteristic of $\N$ is $0$ (this can be checked by substituting $-1$ 
in the known Poincar\'e polynomial of $\N$, see e.g. \cite{atiyah-bott}), 
so the Euler class of $T^*_\N$ vanishes. Therefore $\N$ has self
intersection number $0$ in $T^*_\N\subset \M$. 
This shows that the intersection form is zero.

2. We can also calculate the dimension of the cokernels in our toy example.
Namely, the dimension of $\coker(L_{toy})$ is clearly $1$, as the map 
$L_{toy}:H^0(Z_{toy})\rightarrow H^2(Z_{toy})$ is the multiplication with
$c_1(L_{Z_{toy}})=0$ 
(cf. the example at the end of Section~\ref{kaehler}). 

Thus, by the above theorem, we have that $\coker(j_{\M_{toy}})$ is 
$1$-dimensional. It can be seen directly, using Zariski's lemma 
(Lemma 8.2 in \cite{barth-peters-van} p. 90), that the kernel of the
map $j_{\M_{toy}}$ is generated by the cohomology class of the 
elliptic curve $P$, the generic fibre of the toy Hitchin map, hence it is
$1$ dimensional, indeed.         
 
\end{example}

\end{document}